\documentclass{article}
\usepackage[a4paper,top=3cm,bottom=2cm,left=3cm,right=3cm,marginparwidth=1.75cm]{geometry}

\usepackage{amsmath}
\usepackage{amssymb}
\usepackage{bm}
\usepackage{multirow}
\usepackage{multicol,xcolor}
\usepackage{mathtools}
\usepackage{amsthm}
\usepackage{graphicx}
\usepackage[colorlinks=true, allcolors=blue]{hyperref}
\usepackage{marginnote}
\usepackage{soul}

\usepackage{pgf,tikz}
\usepackage{mathrsfs}
\usetikzlibrary{arrows}

\usepackage[english]{babel}
\usepackage[utf8x]{inputenc}
\usepackage[T1]{fontenc}
\usepackage{enumitem}
\usepackage{bigstrut}
\usepackage{makecell}
\usepackage{bookmark, booktabs}
\usepackage{blkarray}
\usepackage{textcomp}
\usepackage{xcolor}
\usepackage{enumitem}
\usepackage{multicol}
\usepackage{multirow}
\usepackage{bm}
\usepackage{marginnote, mathtools, wasysym}

\usetikzlibrary{decorations.pathreplacing}

\theoremstyle{definition}
\newtheorem{definition}{Definition}[section]
\newtheorem{remark}[definition]{Remark}
\newtheorem{example}[definition]{Example}

\theoremstyle{plain}
\newtheorem{theorem}[definition]{Theorem}
\newtheorem{lemma}[definition]{Lemma}

\newtheorem{property}[definition]{Property}
\newtheorem{conjecture}[definition]{Conjecture}

\renewenvironment{proof}[1][\noindent Proof]{{\par\pushQED{\qed}\itshape #1\@. }}{\popQED}

\DeclareMathOperator{\PG}{PG}
\DeclareMathOperator{\PGL}{PGL}
\DeclareMathOperator{\im}{Im}
\renewcommand{\L}{\mathcal{L}}
\renewcommand{\P}{\mathcal{P}}

\newcommand{\N}{\mathbb{N}}
\newcommand{\F}{\mathbb{F}}

\newcommand{\comments}[1]{}

\newcommand{\qbin}[2]{\genfrac{[}{]}{0pt}{}{#1}{#2}_{q}}
\newcommand{\qbinz}[2]{\genfrac{[}{]}{0pt}{}{#1}{#2}}
\allowdisplaybreaks

\title{Cameron-Liebler sets of generators in the Klein quadric $Q^+(5,q)$}
\author{
 Jozefien D'haeseleer\footnote{Department of Mathematics: Analysis, Logic and Discrete Mathematics, Ghent University, Krijgslaan 281, Building S8, Ghent, Belgium, \url{jozefien.dhaeseleer@ugent.be}}, Jonathan Mannaert\footnote{Department of Mathematics and Data Science, Vrije Universiteit Brussel (VUB), Belgium, \url{Jonathan.Mannaert@vub.be}},
 Leo Storme\footnote{Department of Mathematics: Analysis, Logic and Discrete Mathematics, Ghent University, Krijgslaan 281, Building S8, Ghent, Belgium, \url{leo.storme@ugent.be}}}
\date{}

\begin{document}

\maketitle

\begin{abstract}We investigate Cameron-Liebler sets of planes in the Klein quadric $Q^+(5,q)$ in PG$(5,q)$. We prove that there are many examples of such Cameron-Liebler sets of planes in the Klein quadric. {More specifically, we provide an incomplete list of examples of such Cameron-Liebler sets of planes. By doing so, we also provide some characteristic results regarding these sets in connection with the Klein quadric. These results contribute to an open conjecture posed in \cite{Ferdinand.}}
\end{abstract}

\section{Introduction}\label{sec:intro}
In 1982, Cameron and Liebler introduced particular line sets in $\PG(n,q)$ when investigating the orbits of the action of $\PGL(n+1,q)$ on the points and lines of the projective space $\PG(n,q)$ (\cite{cameronliebler}). {While their main results were obtained in $\PG(3,q)$ it provided the fundamentals for the general case.}
In recognition of their work, these line sets later became known as Cameron-Liebler line sets. 
Over the years, many results concerning the characterisation and classification of Cameron-Liebler sets have been found (see~\cite{bd, cossidentepavese, debeulemannaert, TF:14, gmp, gmet, gmol, met, met2, rod, ThesisDrudge, D99}, among others). 

The many equivalent ways to describe a Cameron-Liebler set, both combinatorially and algebraically, sparked the interest of many researchers, and led to generalisations.
Cameron-Liebler sets of $k$-spaces in $\PG(n,q)$ were studied in~\cite{bdd,bddc,debeulemannaertstorme, Ferdinand.,met3,rodstovan},
and Cameron-Liebler sets of subsets of finite sets were discussed and classified in \cite{meysets}. 
In~\cite{CLpolequiv,CLpol}, (degree one) Cameron-Liebler sets of generators in finite classical polar spaces were introduced, and new examples, using regular ovoids were given in~\cite{CLpolMMJ}.
{Furthermore, Cameron-Liebler sets were also studied in the context of affine-classical spaces~\cite{CLaffineclassical}.}
From a graph theory point of view, Cameron-Liebler sets can be considered as intriguing sets {or 2-perfect colorings} in the sense of~\cite{debruynsuzuki, Ferdinand.}.
The finite set, projective geometry and polar space context correspond to the Johnson, Grassmann and dual polar graphs.
In many contexts, Cameron-Liebler sets correspond to \emph{Boolean degree one functions}, 
{ or completely regular codes of strength $0$, see \cite{Ferdinand.} and the references therein.
 The most common equivalent approach are the Boolean degree one functions of the schemes corresponding to these graph.
 }
\par In each of the contexts, equivalent characterisations of Cameron-Liebler sets exist.
The central question is to classify the Cameron-Liebler sets.
{However, in most context, this has shown to be a complicated problem. For example, in the projective case, plenty of non-existence results exist, but they do not fully solve the problem.
The only non-trivial examples that are found are examples of line classes in $\PG(3,q)$, making this case even more intriguing. Concerning Cameron-Liebler sets in polar spaces, we have a different story. While non-existence results are scarce, the search of general non-trivial examples has been more successful. This provides two non-trivial examples,}
see \cite{CLpolMMJ} and the next section. {Note, that we have to take into account that for polar spaces the problem of Cameron-Liebler sets was approached in two different directions, which coincide in some polar spaces such as $Q^+(5,q)$.} { In \cite{Ferdinand.}, the following conjecture was also posed.
\begin{conjecture}[{\cite[Conjecture 5.1]{Ferdinand.}}]
There exists a constant $n(q,k)$ such that all (degree 1) Cameron-Liebler sets of generators in finite classical polar spaces are the union of the following examples: 
\begin{enumerate}
    \item point-pencils,
    \item all generators in a non-degenerate hyperplane,
    \item all generators in a non-degenerate hyperplane not containing a fixed point in this hyperplane.
\end{enumerate}
\end{conjecture}}

\par In this article, we will focus on new classification results and constructions of non-trivial Cameron-Liebler sets of generators in the hyperbolic quadric $Q^+(5,q)$.  {This work will contribute to the answer of the previous conjecture.}
Section~\ref{sec:prelim} deals with the necessary background of polar spaces, and Cameron-Liebler sets.
In Section \ref{sec:CLunderKleinCorresp}, we discuss Cameron-Liebler sets under the Klein Correspondence, and we use this to give constructions for Cameron-Liebler sets in the Klein quadric, see Subsections \ref{sec:constructionMaxPartSpread}, \ref{sec:BaerSubgeometry}, \ref{sec:LinearSets}. In Section \ref{sec:CharacterisationResultssmallx}, we give a classification result for small values of the parameter $x$, and in the last section, we end with a link between Cameron-Liebler sets and a set of holes of a maximal partial line spread in $\PG(3,q)$.

\section{Preliminaries}\label{sec:prelim}

In this section, we will first {introduce} polar spaces, including the Klein correspondence, and then discuss some background on Cameron-Liebler sets for these polar spaces.

\subsection{Polar spaces}
Polar spaces of rank $d\geq2$ are incidence geometries, whose axiomatic definition goes back to Tits and Veldkamp~\cite{tits,Veldkamp}.

\begin{definition}\label{def:combinatorialpolarspaces}
 A \emph{polar space of rank $d$}, $d\geq3$, is an incidence geometry $(\Pi,\Omega)$ with $\Pi$ a set whose elements are called points and $\Omega$ a set of subsets of $\Pi$ satisfying the following axioms.
 \begin{enumerate}
  \item Any element $\omega\in\Omega$ together with the elements of $\Omega$ that are contained in $\omega$, is a projective geometry of (algebraic) dimension at most $d$.
  \item The intersection of two elements of $\Omega$ is an element of $\Omega$ (the set $\Omega$ is closed under intersections).
  \item For a point $P\in\Pi$ and an element $\omega\in\Omega$ of algebraic dimension $d$ such that $P$ is not contained in $\omega$, there is a unique element $\omega'\in\Omega$ of algebraic dimension $d$ containing $P$ such that $\omega\cap\omega'$ is a hyperplane of $\omega$.
        The element $\omega'$ is the union of all 2-dimensional elements of $\Omega$ that contain $P$ and which intersect $\omega$ in a space of algebraic dimension $1$.
  \item There exist two elements in $\Omega$ both of dimension $d$ whose intersection is empty.
 \end{enumerate}
\end{definition}

We now introduce the finite classical polar spaces.

\begin{definition}
 A \emph{finite classical polar space} is an incidence geometry consisting of the totally isotropic subspaces of a non-degenerate quadratic or non-degenerate reflexive sesquilinear form on a vector space $\F^{n}_{q}$.
\end{definition}

In this article we will consider the finite classical  polar spaces as substructures of a projective space, in which they can naturally be embedded.
We will always use the \emph{projective dimension}.
The subspaces of dimension $0$, $1$ and $2$ are called \emph{points}, \emph{lines} and \emph{planes}, respectively. We now list the 5 finite classical polar spaces of rank $d$.

\begin{itemize}
    \item the hyperbolic quadrics $Q^+(2d-1,q)$ in PG$(2d-1,q)$, 
    \item the elliptic quadrics $Q^-(2d+1,q)$ in PG$(2d+1,q)$,
    \item the parabolic quadrics $Q(2d,q)$ in PG$(2d,q)$,
    \item the Hermitian varieties $H(n,q^2)$ in PG$(n,q^2)$, where $d=\frac{n}{2}$ if $n$ is even and $d=\frac{n+1}{2}$ if $n$ is odd,
    \item the symplectic polar spaces $W(2d-1,q)$ in PG$(2d-1,q)$.
\end{itemize}

\begin{definition}
 The subspaces of maximal dimension (being $d-1$) of a polar space of rank $d$ are called \emph{generators}. We define the \emph{parameter} $e$ of a polar space $\mathcal{P}$ over $\F_{q}$ as $\log_{q}(x-1)$ with $x$ the number of generators through a $(d-2)$-space of $\mathcal{P}$.
\end{definition}

The parameter $e$ of a polar space only depends on the type of the polar space and not on its rank. Table~\ref{parameter} gives an overview.

\begin{table}[ht]
 \centering
 \begin{tabular}{l | c}
  polar space               & $e$ \\ \midrule
  ${Q}^{+}(2d-1,q)$ & 0   \\
  ${H}(2d-1,q)$     & 1/2 \\
  ${W}(2d-1,q)$     & 1   \\
  ${Q}(2d,q)$       & 1   \\
  ${H}(2d,q)$       & 3/2 \\
  ${Q}^{-}(2d+1,q)$ & 2   \\ \bottomrule
 \end{tabular}
 \caption{The parameter $e$ of the finite classical polar spaces.}\label{parameter}
\end{table}

We now present some important facts about the hyperbolic quadric $Q^+(5,q)$ in PG$(5,q)$.

\subsubsection{The hyperbolic quadric \texorpdfstring{$Q^+(5,q)$}{Q+(5,q)} and the Klein correspondence}
A hyperbolic quadric $Q^+(5,q)$ in PG$(5,q)$ is a non-singular quadric with standard equation $X_0X_1+X_2X_3+X_4X_5=0$. It contains points, lines and planes. 

The generators, so planes, of the hyperbolic quadric $Q^+(5,q)$ can be partitioned into two classes, often called the class of the {\em Latin} generators and the class of the {\em Greek} generators. Two generators $\Pi_1$ and $\Pi_2$ of the hyperbolic quadric $Q^+(5,q)$ are equivalent if and only if they are equal or intersect in a point.

There is one particular hyperbolic quadric $Q^+(5,q)$, which plays a special role. This hyperbolic quadric is called the {\em Klein quadric}, because of its relevance for the Klein correspondence.

The lines of the projective space PG$(3,q)$ have Pl\"ucker coordinates $(p_{01},p_{02},p_{03},p_{12},p_{31},p_{23})$. A 6-tuple $(p_{01},p_{02},p_{03},p_{12},p_{31},p_{23})$ is a Pl\"ucker coordinate for a line of PG$(3,q)$ if and only if this 6-tuple satisfies the quadratic condition $p_{01}p_{23}+p_{02}p_{31}+p_{03}p_{12}=0$. Geometrically, this means that a 6-tuple  $(p_{01},p_{02},p_{03},p_{12},p_{31},p_{23})$ is a Pl\"ucker coordinate for a line of PG$(3,q)$ if and only if this 6-tuple defines a projective point of PG$(5,q)$ belonging to the hyperbolic quadric with equation $X_0X_5+X_1X_4+X_2X_3=0$. This particular hyperbolic quadric is called the {\em Klein quadric} of PG$(5,q)$.  The bijective relation between the set of lines of PG$(3,q)$ and the set of points of the Klein quadric, which are the Pl\"ucker coordinates of these lines, is called the {\em Klein correspondence}.

In Table \ref{tabel kleinkwadriek}, we give an overview of the most important correspondences.
\begin{table}[h]
	\centering
	\begin{tabular}{ p{7cm}|p{7cm} }
		
		\boldmath{$\PG(3,q)$} & \boldmath{$Q^+(5,q)$}  \\ \midrule 
		Line  & Point   \\  \midrule
		Two intersecting lines & Two points, contained in a common line of $Q^+(5,q)$ \\ \midrule
		The set of lines through a fixed point $P$ and in a fixed plane $\pi$, with  $P\in \pi$ & Line  \\ \midrule
		The set of lines in a fixed plane & Greek plane 
		\\ \midrule
	The set of lines through a fixed point & Latin plane 
	\\ \bottomrule
	\end{tabular}
 \caption{The image of sets of lines of $\PG(3,q)$ under the Klein correspondence.} \label{tabel kleinkwadriek} 
\end{table}


Consider a point $P$ of the Klein quadric $Q^+(5,q)$. Then to this point $P$, there corresponds a unique tangent hyperplane $T_P(Q^+(5,q))$, which intersects the Klein quadric $Q^+(5,q)$  in a cone with vertex the point $P$ and with base a non-singular hyperbolic quadric $Q^+(3,q)$.

This point $P$ is the Pl\"ucker coordinate of a line $\ell$ of PG$(3,q)$. The $q+1$ points $P_1,\ldots,P_{q+1}$ of this line $\ell$ define the $q+1$ Latin generators through this point $P$, and the $q+1$ planes $\pi_1,\ldots,\pi_{q+1}$ of PG$(3,q)$ through this line $\ell$ define the $q+1$ Greek generators through this point $P$. Together, these Latin and Greek generators define the set of $2(q+1)$ generators of the Klein quadric through this point $P$. We call this set of the $2(q+1)$ generators through a point $P$ of the Klein quadric, the {\em point-pencil} with vertex $P$.\\

We now introduce the \emph{Gaussian binomial coefficient}, which is very useful to describe counting results for vector spaces and polar spaces.

\begin{definition}
 The \emph{Gaussian binomial coefficient} for integers $n,k$, with $n\geq k\geq0$ and prime power $q\geq2$, is given by
 \[
  \qbin{n}{k}=\prod^{k}_{i=1}\frac{q^{n-k+i}-1}{q^{i}-1}=\frac{(q^{n}-1)\cdots(q^{n-k+1}-1)}{(q^{k}-1)\cdots(q-1)}\:.
 \]
 For integers $n,k$, with $k>n\geq0$ or $n\geq0>k$, we set $\qbin{n}{k}=0$.
\end{definition}

The Gaussian binomial coefficient $\qbin{n}{k}$ equals the number of $(k-1)$-spaces in the projective space $\PG(n-1,q)$.
Moreover, we will denote the number $\qbin{n+1}{1}$
of points in PG$(n, q)$ by the symbol $\theta_n(q)$. 
If the size $q$ of the field $\mathbb{F}_q$ is clear from the context, we will write $\qbin{n}{k}$ as $\qbinz{n}{k}$.

\subsection{Cameron-Liebler sets for polar spaces}\label{sectionCLintro}
{In general, two objects, that often coincide, are known as Cameron-Liebler sets of generators. The first object refers to the algebraic underlying structure and equals the definition of Boolean degree $1$ functions in these polar spaces. These sets are what we call degree one Cameron-Liebler sets and are described in} \cite{CLpolequiv, CLpol}. {In order to define these sets, we need the notion of the point-generator incidence matrix $A$. For this $\{0,1\}$-valued matrix, the rows and columns of this matrix correspond to the points and generators, respectively. A certain position $(i,j)$ contains 1 if and only if the corresponding point is incident with the corresponding generator.}

\begin{definition}\label{defspecialCL}
    {Let $A$ be the point-generator incidence matrix.} A degree one Cameron-Liebler set of generators  in a finite classical polar space $\mathcal{P}$ is a set {$\mathcal{L}$} of generators in $\mathcal{P}$, with characteristic vector $\chi$ such that $\chi \in \im(A^T)$. {Moreover, we say that $\mathcal{L}$ has parameter $x$ if $$|\mathcal{L}|=x\prod_{i=0}^{d-2}(q^{e+i}+1),$$
    where $d$ denotes the rank of the polar space and $e$ the type.} 
\end{definition}

{The close connection with the algebraic nature of these objects has certain advantages. Moreover, regarding the underlying dual polar scheme, it is shown that these degree $1$ functions/sets  have characteristic vectors contained in the first two eigenspaces $V_0\perp V_1$, see \cite{CLpolMMJ, Ferdinand.}. For more information on association schemes, we refer the reader to \cite{BCN}.}
{The second object is defined using the geometrical approach. This definition lies in the perspective of equitable bipartitions or intruiging sets, For more information regarding this approach we refer the reader to \cite[Section 1.3]{DBDIM23}.  This definition is based on the \emph{disjointness-definition}.}

\begin{definition}[{\cite{CLpol}}]\label{defCL}
	Let $\mathcal{P}$ be a finite classical polar space with parameter $e$ and rank $d$. A set $\mathcal{L}$ of generators in $\mathcal{P}$ is a Cameron-Liebler set of generators in $\mathcal{P}$, with parameter $x$, if and only if for every generator $\pi$ in $\mathcal{P}$, the number of elements of $\mathcal{L}$, disjoint from $\pi$, equals $(x-\chi(\pi))q^{\binom{d-1}{2}+e(d-1)}$. 
\end{definition}
Using association scheme notation, we can interpret the previous definition as follows. The characteristic vector of a Cameron-Liebler set is contained in $V_0 \perp W$, with $W$ the eigenspace of the disjointness matrix $A_d$ corresponding to a specific eigenvalue. It can be seen that $W$ always contains $V_{1}$, but it does not necessarily coincide with $V_{1}$. Clearly, each degree one Cameron-Liebler set is a Cameron-Liebler set, but not vice versa.

Moreover, for most polar spaces the definition of a Cameron-Liebler set and a degree one Cameron-Liebler set coincide; the exceptions are the hyperbolic quadrics of even rank, the parabolic quadrics of odd rank and the symplectic polar spaces of odd rank, i.e.~the polar spaces $Q^{+}(4n-1,q)$, $Q(4n+2,q)$ and $W{+}(4n+1,q)$ for all $q$ and $n$. 
See~\cite[Section 3]{CLpol} for a detailed discussion.

Hence, note that for this paper, in which we work in $Q^+(5,q)$, we do not have to distinguish between the two types of Cameron-Liebler sets, as they are the same.

 We also have the following result, which will be very useful in this paper.

\begin{theorem}[{\cite[Theorem 3.1]{CLpolequiv}}]\label{stellinggg}
 Let $\mathcal{P}$ be a 
 polar space of rank $d$ with parameter $e$ defined over $\F_{q}$, let $\mathcal{L}$ be a set of generators of $\mathcal{P}$ and let $i$ be an integer with $1\leq i\leq d$.
 If $\mathcal{L}$ is a degree one Cameron-Liebler set of generators in $\mathcal{P}$, with parameter $x$,  then the number of elements of $\mathcal{L}$ meeting a generator $\pi$ in a $(d-i-1)$-space equals
 \begin{align}\label{formulelang}
  \begin{dcases}
   \left( (x-1)\begin{bmatrix} d-1 \\ i-1 \end{bmatrix} + q^{i+e-1}\begin{bmatrix} d-1 \\ i \end{bmatrix}\right)q^{\binom{i-1}{2}+ (i-1)e}
    & \text{if $\pi \in \mathcal{L}$}     \\
   \hfil x \begin{bmatrix} d-1 \\ i-1\end{bmatrix} q^{\binom{i-1}{2}+(i-1)e}
    & \text{if $\pi \notin \mathcal{L}$.}
  \end{dcases}
 \end{align}
 Moreover, if this property holds for a polar space $\mathcal{P}$ and an integer $i$ such that
 \begin{itemize}
  \item $i$ is odd if $\mathcal{P}=Q^+(2d-1,q)$,
  \item $i\neq d$ if $d$ is odd, and $\mathcal{P}=Q(2d,q)$ or $\mathcal{P}=W(2d-1,q)$,
  \item $i$ is arbitrary if $\mathcal{P}$ is another polar space,
 \end{itemize}
 then $\mathcal{L}$ is a degree one Cameron-Liebler set with parameter $x$.
\end{theorem}

For $\P = Q^+(5,q)$, we have that $e=0$ and $d=3$. Hence, we have the following intersection numbers.
\begin{align*}
	\left\{ \begin{matrix}
	\left( (x-1) \begin{bmatrix}2 \\i-1\end{bmatrix} +q^{i-1}\begin{bmatrix}2 \\ i\end{bmatrix} \right) q^{\binom{i-1}{2}} & \mbox{if } \pi \in \mathcal{L}\\  x \begin{bmatrix} 2 \\ i-1\end{bmatrix} q^{\binom{i-1}{2}} & \mbox{if }\pi \notin \mathcal{L}. 
	\end{matrix}\right.
	\end{align*}
This is summarized in Table \ref{tabel intersection}.
\begin{table}[ht]
 \centering
\begin{tabular}{c | c c c} 

 Distance & $i=1$ & $i=2$ & $i=3$ \\ [0.5ex] 
 \midrule
 $\pi\in \L$ & $x+q$ & $xq+x-1$ & $(x-1)q$ \\ 

 $\pi\notin \L$ & $x$ & $x(q+1)$ & $x$ \\ [1ex] 
 \bottomrule
\end{tabular}\caption{Intersection numbers planes in $Q^+(5,q)$} \label{tabel intersection} 
\end{table}

We also give some properties of degree one Cameron-Liebler sets of generators in polar spaces that can easily be proven. 

\begin{lemma}[{\cite[Lemma 3.3]{CLpolequiv} and~\cite[Lemma 4.1]{CLpol}}]\label{basicoperations}
 Let $\mathcal{P}$ be a 
 polar space of rank $d$, and with the set of generators $\Omega_d$. Let $\mathcal{L}$ and $\mathcal{L}'$ be (degree one) Cameron-Liebler sets of $\mathcal{P}$ with parameters $x$ and $x'$ respectively.
 \begin{itemize}
  \item[(i)] $0\leq x\leq q^{e+d-1}+1$.
  \item[(ii)] The complement of $\mathcal{L}$ in $\Omega_{d}$ is a (degree one) Cameron-Liebler set of $\mathcal{P}$ with parameter $q^{e+d-1}+1-x$.
  \item[(iii)] If $\mathcal{L}\cap\mathcal{L}'=\emptyset$, then $\mathcal{L}\cup\mathcal{L}'$ is a (degree one) Cameron-Liebler set of $\mathcal{P}$ with parameter $x+x'$.
  \item[(iv)] If $\mathcal{L}'\subset\mathcal{L}$, then $\mathcal{L}\setminus\mathcal{L}'$ is a (degree one) Cameron-Liebler set of $\mathcal{P}$ with parameter $x-x'$.
 \end{itemize}
\end{lemma}

We now present two examples of (degree one) Cameron-Liebler sets.

\begin{example}[{\cite[Example 4.2]{CLpol}}]\label{ex:pointpencil}
 A \emph{point-pencil} with \emph{vertex} $P$ in a polar space $\mathcal{P}$ is the set of all generators in $\mathcal{P}$ containing the point $P$. A point-pencil is a degree one Cameron-Liebler set (and thus also a Cameron-Liebler set) since { the characteristic vector is a column of $A$.} 
 It has parameter 1.
 \par {A second example can be found in \cite{CLpol} and uses partial ovoids in the construction.} A \emph{partial ovoid} is a set of points on $\mathcal{P}$ which are pairwise not collinear. In other words, for any two points the line through them is not isotropic with respect to the underlying quadratic or sesquilinear form.
 Hence, for any two points of a partial ovoid, the point-pencils with these vertices are disjoint. So, if $\mathcal{O}$ is a partial ovoid, the union of the point-pencils with the points of $\mathcal{O}$ as vertex, is a (degree one) Cameron-Liebler set with parameter $|\mathcal{O}|$.
\end{example}

\begin{example}[{\cite[Example 4.4]{CLpol}}]\label{ex:embedded}
 {Let $\mathcal{P}$ and $\mathcal{P}'$ both be polar spaces of rank $d$ and suppose that $\mathcal{P}$ has parameter $e$ and $\mathcal{P}'$ has parameter $e-1$ and is embedded in $\mathcal{P}$. Only three such examples exist, see the list below. Then $\mathcal{P}'$ is a}
 Cameron-Liebler set with parameter $q^{e-1}+1$.
 \par 
 {A more specific list (up to a projection) of such examples is described below:}
 \begin{itemize}
  \item There are parabolic quadrics $Q(2d,q)$ embedded in the elliptic quadric $\mathcal{P}=Q^{-}(2d+1,q)$. Each of them gives rise to a (degree one) Cameron-Liebler set of $\mathcal{P}$ with parameter $q+1$.
  \item There are hyperbolic quadrics $Q^{+}(2d-1,q)$ embedded in the parabolic quadric $\mathcal{P}=Q(2d,q)$. Each of them gives rise to a (degree one) Cameron-Liebler set of $\mathcal{P}$ with parameter $2$.
        Recall that the symplectic variety $W(2d-1,q)$ is isomorphic to $\mathcal{P}$ if $q$ is even.
  \item There are Hermitian polar spaces $H(2d-1,q)$ embedded in the Hermitian polar space $\mathcal{P}=H(2d,q)$, $q$ a square.
        Each of them gives rise to a (degree one) Cameron-Liebler set of $\mathcal{P}$ with parameter $\sqrt{q}+1$.
 \end{itemize}
\end{example}

The (degree one) Cameron-Liebler sets from Examples~\ref{ex:pointpencil} and~\ref{ex:embedded}, their complements and disjoint unions, are called \emph{trivial}; all others are called non-trivial.

Recently, the first non-trivial degree one Cameron-Liebler sets were found. In \cite{CLpolMMJ}, two constructions are given:  they first give a construction of a family of non-trivial Cameron-Liebler sets of generators in the hyperbolic quadrics $Q^+(5, q)$, $q$ odd, in which they also use the Klein correspondence with the lines in $\PG(3,q)$. 
Furthermore, they also introduce
and discuss a new construction of Cameron-Liebler sets for polar spaces with $e \leq 1$ and rank $d\geq 4$. The latter is based on a generalization of ovoids to higher-dimensional subspaces of polar spaces.\\

\subsection{{Some known classification results}}
Finally we mention a few classification results. Note that from Definition~\ref{defCL} and Theorem \ref{stellinggg}, {the following result can be deduced.} 

\begin{theorem}[{\cite[Lemma 4.8]{CLpol}} and {\cite[Lemma 5.3]{CLpolequiv}}]
 Let $\mathcal{P}$ be a 
 polar space. If $\mathcal{L}$ is a (degree one) Cameron-Liebler set of $\mathcal{P}$ with parameter $x$, then $x\in\N$.
\end{theorem}

The most important classification result so far states that small (degree one) Cameron-Liebler sets are trivial.

\begin{theorem}[{\cite[Theorem 5.5]{CLpolequiv}}]\label{th:classificationsmalldegreeoneCL}
 Let $\mathcal{P}$ be a 
 polar space of rank $d$ and parameter $e$, and let $\mathcal{L}$ be a degree one Cameron-Liebler set of $\mathcal{P}$ with parameter $x$.
 If $x\leq q^{e-1}+1$, then $\mathcal{L}$ is the union of $x$ point-pencils whose vertices are pairwise non-collinear or $x=q^{e-1}+1$ and  $\mathcal{L}$ is the set of generators in an embedded polar space of rank $d$ and with parameter $e-1$.
\end{theorem}

\begin{theorem}[{\cite[Theorem 6.7]{CLpol}}]\label{th:classificationsmallCL}
 Let $\mathcal{P}$ be a 
 polar space of rank $d$ and parameter $e$, which is not a hyperbolic quadric of even rank, a parabolic quadric of odd rank or a symplectic polar space of odd rank.
 Let $\mathcal{L}$ be a Cameron-Liebler set of $\mathcal{P}$ with parameter $x$. If $x\leq q^{e-1}+1$, then $\mathcal{L}$ is the union of $x$ point-pencils whose vertices are pairwise non-collinear or $x=q^{e-1}+1$ and $\mathcal{L}$ is the set of generators in an embedded polar space of rank $d$ and with parameter $e-1$.
\end{theorem}

From Theorem~\ref{th:classificationsmalldegreeoneCL}, it follows that all degree one Cameron-Liebler sets with parameter 1 are point-pencils. {Note that this theorem is only valid for degree $1$ Cameron-Liebler sets and not for the geometric equivalent described by Definition \ref{defCL}.}
In~\cite[Theorem 6.4]{CLpol}, Cameron-Liebler sets with parameter 1 that are not trivial are described, showing that the exclusion of several types of polar spaces in Theorem~\ref{th:classificationsmallCL} is necessary.
\par For the symplectic polar space $W(5,q)$ and the parabolic quadric $Q(6,q)$, it is proven in~\cite[Theorem 5.9]{CLpolequiv} that a Cameron-Liebler set with parameter $x$, with $2\leq x\leq \sqrt[3]{2q^2}-\frac{\sqrt[3]{4q}}{3}+\frac{1}{6}$, is a union of embedded polar spaces $Q^+(5,q)$ and point-pencils.

\begin{remark}
{Cameron-Liebler sets of generators and (degree $1$) Cameron-Liebler sets of generators are equal objects in $Q^+(5,q)$, see Section \ref{sectionCLintro} and  ~\cite[Section 3]{CLpol} for a detailed discussion. For the main results on degree $1$ Cameron-Liebler sets in $Q^+(5,q)$,}
    we see that the only known classification results give 
    that the parameter $x\in \mathbb{N}$, and that for $x=1$,     
    the corresponding Cameron-Liebler set is a point-pencil. 
    Furthermore, only two non-trivial Cameron-Liebler sets are known, see \cite{CLpolMMJ}.
\end{remark}


In this paper, we show that there are many examples of Cameron-Liebler sets of planes in $Q^+(5,q)$, and we also present some characterisation results.

\section{Cameron-Liebler sets under the Klein correspondence}\label{sec:CLunderKleinCorresp}

Let $\L$ be a Cameron-Liebler set of planes in $Q^+(5,q)$ with parameter $x$.  {Then $|\L| = 2x(q+1)$, and, from Theorem \ref{stellinggg}, it follows that given a Greek (resp. Latin) plane $\pi\in \L$, there are $x(q+1)-1$ Greek (resp. Latin) planes in $\L$ meeting $\pi$ in precisely one point. Hence, $\L$ consists of $x(q+1)$ Greek and $x(q+1)$ Latin planes.}
This set $\L$ defines, using the Klein correspondence, in $\PG(3,q)$ a union $P_0\cup P_2$ of two sets in which $P_0$ is a set of $x(q+1)$ points, and $P_2$ is a set of $x(q+1)$ planes, with the following property.
\begin{property}\label{prop}
{Suppose that $P_0$ is a set of points and $P_2$ is a set of planes in $\PG(3,q)$, respectively. Then the following statements are equivalent:}
\begin{enumerate}
    \item The set of generators in $Q^+(5,q)$ derived from $P_0$ and $P_2$, using the Klein correspondence is a degree $1$ Cameron-Liebler set of parameter $x$.
    \item The following two properties are valid.
    \begin{itemize}
    \item Every plane of $\PG(3,q)$ contains $x$ or $q+x$ points of $P_0$, and the planes of $P_2$ are the planes containing $q+x$ points of $P_0$.
    \item Every point of $\PG(3,q)$ lies in $x$ or $q+x$ planes of $P_2$, and the points of $P_0$ are the points lying in $q+x$ planes of $P_2$.
\end{itemize}
\end{enumerate}

\end{property}
\begin{proof}
    Each plane and each point in $\PG(3,q)$ correspond to a Greek and Latin plane in $Q^+(5,q)$ respectively. An incident point-plane pair in $\PG(3,q)$ corresponds to two planes in $Q^+(5,q)$ at distance 1. Then the properties follow from Theorem \ref{stellinggg}, see the first column in Table \ref{tabel intersection}. {In fact, from Theorem \ref{stellinggg}, we know that these properties form a sufficient condition for degree $1$ Cameron-Liebler sets of parameter $x$.}\qedhere \\
\end{proof}


\subsection{Cameron-Liebler sets arising from partial line spreads in  \texorpdfstring{$\PG(3,q)$}{PG(3,q)}}\label{sec:constructionMaxPartSpread}



The first class of examples,  giving already a  Cameron-Liebler set of generators of the Klein quadric, for every parameter $x$, with $1\leq x \leq q^2+1$, is obtained via the point-pencils of pairwise non-collinear vertices. We present the construction in PG$(3,q)$, and then use the Klein correspondence to define the corresponding Cameron-Liebler set on the Klein quadric.
\begin{example}
    Take $x$ pairwise disjoint lines $\ell_1, \ldots, \ell_x$ in $\PG(3,q)$. Take the set of the $x(q+1)$ points on these lines, together with the set of $x(q+1)$ planes through these lines. This gives a {Cameron-Liebler} set with parameter $x$, which corresponds to the union of $x$ point-pencils in $Q^+(5,q)$, whose points are pairwise non-collinear.
\end{example}



This example is linked to the theory of partial spreads and spreads in $\PG(3,q)$.

\begin{definition}
A set of $x$ pairwise disjoint lines in $\PG(3,q)$ is called a {\em partial spread of size $x$}.

A partial spread in $\PG(3,q)$ is called {\em maximal} when it is not contained in a larger partial spread of $\PG(3,q)$.

A {\em spread} of $\PG(3,q)$ is a set of $q^2+1$ pairwise disjoint lines, which then form a partition of the point set of $\PG(3,q)$.

A {\em hole} of a partial spread of $\PG(3,q)$ is a point, not belonging to a line of the partial spread.
\end{definition}

There exist many spreads in $\PG(3,q)$. The classical example is the {\em regular} spread of $\PG(3,q)$. The following spectrum results on maximal partial spreads in $\PG(3,q)$ are known.

\begin{remark}
    Maximal partial spreads $S$ exist in $\PG(3,q)$ for 
    \begin{itemize}
        \item $q$ odd, $q\geq 7$, $\frac{q^2+1}{2}+6\leq |S| \leq q^2-q+2$. 
        \item $q$ even, $q\geq q_0$, $\frac{5q^2+q+16}{8}\leq |S|\leq q^2-q+2$.
    \end{itemize}

    The results for $q$ odd were proven by Heden in \cite{MR1235896,MR1341439,MR1874731} and for $q$ even by Heden and Storme \cite{HedenStorme}.
\end{remark}

\begin{example}\label{ex:dual} Let $S$ be a maximal partial spread of size $q^2+1-x$. The set of holes $P_0$ is of size $x(q+1)$, and {dually} the set  of planes $P_2$, not containing a line of $S$, has size $x(q+1)$. {Using the Klein correspondence, $P_0\cup P_2$ provides a  Cameron-Liebler set of generators in $Q^+(5,q)$. }
\end{example}
\begin{proof}
{Suppose that the assertion is valid. Note that a plane $\pi$ of $\PG(3,q)$ contains $x$ or $q+x$ holes, and dually, a point $P\in \PG(3,q)$ is contained in $x$ or $q+x$ planes not containing a line of the spread. 
By Property \ref{prop}, we find that the set $P_0\cup P_2$ forms a set of points and planes, which gives a Cameron-Liebler set in $Q^+(5,q)$, with parameter $x$, under the Klein correspondence.} 
\end{proof}
\begin{remark}
    Note that Example \ref{ex:dual} is the complement of the Cameron-Liebler set arising from the maximal partial spread $S$. {Note that the important fact here is that this example corresponds with a set of holes of a maximal partial spread. This implies that the example does not contain any point-pencils.}
   This follows as a point-pencil in $Q^+(5,q)$ corresponds to a set of $q+1$ points on a line $\ell$ in $\PG(3,q)$, and the set of planes through this line. If this would be contained in this Cameron-Liebler set, then the line $\ell$ would extend the maximal partial spread, a contradiction.
\end{remark}

\subsection{Cameron-Liebler sets arising from a Baer subgeometry}\label{sec:BaerSubgeometry}

We now describe an example of a Cameron-Liebler set with parameter $q+1$ on the Klein quadric, arising from a Baer subgeometry $\PG(3,q)$ in $\PG(3,q^2)$.

\begin{example}
    The set of points and the set of planes of the Baer subgeometry $\PG(3,q)$ defines a {Cameron-Liebler} set with parameter $q+1$ on the Klein quadric $Q^+(5,q^2)$. In fact, this Cameron-Liebler set with parameter $q+1$ on the Klein quadric corresponds to a Klein sub hyperbolic quadric $Q^+(5,q)$ in the Klein quadric $Q^+(5,q^2)$. We will call this type of Cameron-Liebler set with parameter $q+1$, a Cameron-Liebler set of {\em  Baer subgeometry type}.
\end{example}
\begin{proof}
Consider a Baer subgeometry $\PG(3,q)$, naturally embedded in $\PG(3,q^2)$. Then the planes of $\PG(3,q^2)$ either intersect this Baer subgeometry $\PG(3,q)$  in $q+1$ points or in $q^2+q+1$ points. In the first case, the plane only shares a Baer subline with the Baer subgeometry $\PG(3,q)$, and in the second case, the plane  shares a Baer subplane with the Baer subgeometry $\PG(3,q)$.

Similarly, dually, the points of $\PG(3,q^2)$ lie in $q+1$ planes of the Baer subgeometry $\PG(3,q)$, or in $q^2+q+1$ planes of the Baer subgeometry $\PG(3,q)$. In fact, a point $P$ of the   Baer subgeometry $\PG(3,q)$  lies in $q^2+q+1$ planes of the Baer subgeometry $\PG(3,q)$, while a point $P$ of $\PG(3,q^2)$, not belonging to the Baer subgeometry $\PG(3,q)$, only lies in $q+1$ planes of the Baer subgeometry $\PG(3,q)$. These latter $q+1$ planes of the  Baer subgeometry $\PG(3,q)$  are the $q+1$ planes of the  Baer subgeometry $\PG(3,q)$ through the line $PP^{q}$, defined by the point $P$ and its conjugate point $P^{q}$ with respect to the Baer subgeometry $\PG(3,q)$. {The assertion follows from Property \ref{prop}.}
\end{proof}

\subsection{ Cameron-Liebler sets arising from linear sets}\label{sec:LinearSets}

We now present another example arising from a subgeometry. As we will state later, it is part of a class of Cameron-Liebler sets arising from linear sets. {We will  show the following theorem.}
\begin{theorem}\label{ex:linsets}  There exist Cameron-Liebler sets on the Klein quadric $Q^+(5,q^t)$ with parameter $x=\frac{q^t-1}{q-1}$, arising from scattered $\mathbb{F}_q$-linear sets of rank $\frac{rt}{2}$.
\end{theorem}

    Embed $\PG(5,q)$ in $\PG(5,q^3)$. Take a line $\ell$ in $\PG(5,q^3)$ skew to $\PG(5,q)$, with $\ell^q$ and $\ell^{q^2}$ the conjugate lines of the line $\ell$ with respect to the subgeometry $\PG(5,q)$, and such that $\{\ell, \ell^q,\ell^{q^2}\}$ defines a regular $2$-spread in $\PG(5,q)$, or, equivalently, $\dim\langle \ell,\ell^q,\ell^{q^2}\rangle=5$. Project $\PG(5,q)$ from the line $\ell$ to a $3$-dimensional space $\alpha:=\PG(3,q^3)$ skew to $\ell$.
    
    Take a point $P\in \ell$ and its conjugate points $P^q,P^{q^2}$ with respect to the subgeometry $\PG(5,q)$. Then $\langle P, P^q, P^{q^2}\rangle$ is a plane in $\PG(5,q^3)$, which shares a subplane $\PG(2,q)$ with $\PG(5,q)$.  
    Now, every plane $\langle P, P^q, P^{q^2}\rangle$, defined by $P\in \ell$, lies in $q^2+q+1$ different $3$-spaces $\pi$ sharing $\PG(3, q)$ with $\PG(5, q)$. Such a solid $\pi$ does not intersect $\ell$ in a second point, as otherwise $\ell\subset \pi$, and since $\pi$ is defined over $\mathbb{F}_q$, also $\ell^q$ and $\ell^{q^2}$ would lie in $\pi$, a contradiction, since these three lines span the $5$-space.

    So $\ell$ and such a solid $\pi$ span a $4$-space $\PG(4,q^3)$. This $4$-space is projected from $\ell$ onto a plane of $\alpha =\PG(3,q^3)$, containing $q^3+q^2+q+1$ projected points. 

    So we already have $(q^3+1)(q^2+q+1)$ planes of $\alpha$ sharing $q^3+q^2+q+1$ points with the projected subgeometry $\PG(5,q)$ in $\alpha$.\\

    Vice versa, every 4-space $\Pi_4$ of $\PG(5,q^3)$ through the line $\ell$, which shares a  3-space $\PG(3,q)$ with the  subgeometry $\PG(5,q)$ in $\PG(5,q^3)$, is of this type. For extend this 3-dimensional subgeometry $\PG(3,q)$ to $\PG(3,q^3)$ in $\PG(5,q^3)$. Then this 3-dimensional subspace $\PG(3,q^3)$ shares a point $R$ with $\ell$ since they both lie in $\Pi_4$. If we then apply conjugation with respect to the subgeometry $\PG(5,q)$ in $\PG(5,q^3)$, then this 4-space $\Pi_4$ also contains the conjugate points $R^q$ and $R^{q^2}$, since the subspace $\PG(3,q^3)$ stays invariant under this conjugation. So there is a plane $\langle R,R^q,R^{q^2}\rangle$ contained in this subspace $\PG(3,q^3)$, sharing a subplane $\PG(2,q)$ with the subgeometry $\PG(3,q)$ in $\PG(3,q^3)$. This proves that this 4-space $\Pi_4$ is indeed of the type we have discussed.

    No 4-space $\Pi_4$ of $\PG(5,q^3)$ through the line $\ell$ can intersect the subgeometry $\PG(5,q)$ in a 4-dimensional space $\PG(4,q)$. Because then this 4-space would be invariant under the conjugation with respect to the subgeometry $\PG(5,q)$. Hence, this 4-space would contain the three lines $\ell$, $\ell^q$ and $\ell^{q^2}$. This is false since these three lines
    $\ell, \ell^q$ and $\ell^{q^2}$ generate 5 dimensions.

    This implies that all other hyperplanes $\Pi_4$ through the line $\ell$ intersect the subgeometry in a plane $\pi_2=\PG(2,q)$, which is contained in the plane $\Pi_4\cap \Pi_4^q\cap \Pi_4^{q^2}$. 

    This then implies that if we project the subgeometry $\PG(5,q)$ from the line $\ell$ into this 3-space $\alpha=\PG(3,q^3)$, skew to the line $\ell$, we obtain a projected subgeometry $\PG(5,q)$ of size $(q^3+1)(q^2+q+1)$ in $\PG(3,q^3)$, intersecting every plane in either $q^2+q+1$ or $q^3+q^2+q+1$ points. Standard calculations again show that a point of $\alpha$ either lies in $q^2+q+1$ or $q^3+q^2+q+1$ planes of $\alpha$, containing $q^3+q^2+q+1$ projected points. The points of $\alpha$ belonging to $q^3+q^2+q+1$ planes of $\alpha$, containing $q^3+q^2+q+1$ projected points, are precisely, the projected points of the projected subgeometry $\PG(5,q)$.

    Hence, the set of points $P_0$ of the  projected points of the projected subgeometry $\PG(5,q)$, and the set of planes $P_2$ each containing $q^3+q^2+q+1$ projected points define a Cameron-Liebler set on the Klein quadric $Q^+(5,q^3)$ with parameter $x=q^2+q+1$. We call this type of Cameron-Liebler set of {\em projected $\PG(5,q)$ type}.

    The preceding example is a particular example of the following set of examples of Cameron-Liebler sets on the Klein quadric. This set of examples arises from the theory of linear sets.

    We will, in particular, use the following theorem, to make the link between Cameron-Liebler sets in the Klein quadric, and linear sets. 
    \begin{theorem}[{\cite[Theorem 4.2]{MR1772206}}] Let $W_{rt/2}$ be a subspace of rank $rt/2$ of the vector space $V(rt,q)$ of dimension $rt$ over the finite field of order $q$, which is scattered with respect to a normal $t$-spread $S$ of $V(rt,q)$, then the corresponding linear set $B(W_{rt/2})$ is a $2$-intersection set of $\PG(r-1,q^t)$ with respect to the hyperplanes, and with intersection numbers $\theta_{\frac{rt}{2}-t-1}(q)$
and  $\theta_{\frac{rt}{2}-t}(q)$.   \end{theorem}

We  {want} to apply this result in $\PG(3,q^t)$, so with $r=4$.
Then $W_{rt/2}=W_{2t}$ is of rank $2t$ in $V(4t,q)$, $\theta_{\frac{rt}{2}-t-1}(q)=\theta_{t-1}(q)=\frac{q^t-1}{q-1}$ and  $\theta_{\frac{rt}{2}-t}(q)=\theta_{t}(q)=\frac{q^{t+1}-1}{q-1}$.


Geometrically, this means that a subgeometry $\PG(2t-1,q)$, naturally embedded in $\PG(4t-1,q^t)$, is projected from a suitably chosen subspace $\PG(4t-5,q^t)$ into a 
scattered linear set of this 3-space $\PG(3,q^t)$.

Then the size of this scattered linear set is $|P_0|=|\PG(2t-1,q)|=\frac{q^{2t}-1}{q-1}$ and every plane of $\PG(3,q^t)$ contains either $\theta_{t-1}(q)$ or $\theta_t(q)$ points of the set $P_0$.

There are $y=|\PG(2t-1,q)|$ planes of $\PG(3,q^t)$ containing $\theta_t(q)$ points of the set $P_0$. 

This follows from the following equation in which we count the incidences of the points of $\PG(3,q^t)$ with the planes of $\PG(3,q^t)$:

\[y\theta_t(q)+ (|\PG(3,q^t)|-y)\theta_{t-1}(q)=|\PG(2t-1,q)|(q^{2t}+q^t+1).\]

Let $P_2$ be the set of $|\PG(2t-1,q)|$ planes of $\PG(3,q^t)$
containing $(q^{t+1}-1)/(q-1)$ points of the set $P_0$.

The dual similar calculations show that a point of $P_0$ belongs to $(q^{t+1}-1)/(q-1)$ planes of the set $P_2$, while a point of $\PG(3,q^t)\setminus P_0$ belongs to $(q^t-1)/(q-1)$ planes of the set $P_2$.

This shows that this set of points $P_0$ and this set of planes $P_2$ satisfy Property \ref{prop}, hence they define a Cameron-Liebler set of the Klein quadric $Q^+(5,q^t)$ 
with parameter $x=(q^t-1)/(q-1)$.

The particular example of the scattered projected subgeometry $\PG(5,q)$ in $\PG(3,q^3)$ is such a Cameron-Liebler set with parameter $\frac{q^3-1}{q-1}=q^2+q+1$.

In his PhD thesis, see \cite{PhDLavrauw}, M. Lavrauw proved that such an example exists for all the possible parameters. In \cite[Theorem 2.5.5]{PhDLavrauw}, M. Lavrauw proved that when $r$ is even, there always exist scattered $\mathbb{F}_q$-linear sets of rank $\frac{rt}{2}$ in $\PG(r-1,q^t)$, for $t\geq 2$.
{This then leads to Theorem \ref{ex:linsets} on Cameron-Liebler sets of generators on the Klein quadric.}

\subsection{Some characterisation results for small \texorpdfstring{$x$}{x}}\label{sec:CharacterisationResultssmallx}

We first present a characterisation result for Cameron-Liebler sets on the Klein quadric,  {for small values of the parameter $x$.} 

\begin{theorem} Every Cameron-Liebler set $\L$ on the Klein quadric, with parameter $x$ satisfying $1\leq x<\sqrt{q}+1$, is the union of $x$ point-pencils, defined by $x$ points 
pairwise non-collinear on the Klein quadric.
\end{theorem}

\begin{proof}
We present the proof in $\PG(3,q)$.
Recall that any point $p\in \PG(3,q)\setminus P_0$ is contained in $x$ planes of $P_2$.
Consider a point $p\in P_0$. This point is contained in $q+x$ planes of $P_2$, and each of these planes contains $q+x-1$ points of $P_0$, different from $p$. Moreover, there are $x(q+1)-1$ points of $P_0$, different from $p$. By double counting the set $\{(\pi, p') | \pi\in P_2, p'\in P_0, p\in \pi, p'\in \pi\setminus \{p\} \}$, we find that there is some point $p'\in P_0\setminus \{p\}$, which lies in at least $\frac{(q+x)(q+x-1)}{x(q+1)-1}$ planes of $P_2$ through $p$.

Consider the line $\langle p,p'\rangle$. Suppose that there is a point $p''\in \langle p,p'\rangle$, which is not contained in $P_0$. This point $p''$ lies in $x$ planes of $P_2$, and hence, $\frac{(q+x)(q+x-1)}{x(q+1)-1}\leq x$, which implies that $x^2-2x-q+1\geq 0$. But this contradicts  $x< \sqrt{q}+1$. So all points of the line $\langle p,p'\rangle$ belong to $P_0$. 

This then implies that all the planes of $\PG(3,q)$ through the line $\langle p,p'\rangle$ contain already at least $q+1$ points of $P_0$. Since $x<\sqrt{q}+1$, these planes then necessarily contain $q+x$ points of the set $P_0$. Hence, they all belong to the set $P_2$.

This shows that the Cameron-Liebler set  defined by the point-pencil of the Pl\"ucker coordinate of the line $\langle p,p'\rangle$ is contained in the Cameron-Liebler set $\L$. This point-pencil is a Cameron-Liebler set with parameter one that can be removed from $\L$, see Lemma \ref{basicoperations}$(iv)$. This then gives a new Cameron-Liebler set with parameter $x-1$.

Repeating this argument, this shows that the Cameron-Liebler set $\L$ is the union of $x$ pairwise disjoint point-pencils. Since these $x$ point-pencils are pairwise disjoint, their vertices 
are necessarily pairwise non-collinear on the Klein quadric.
\end{proof}


\begin{remark}
    We can use a similar argument to find properties of Cameron-Liebler sets in other polar spaces.
\end{remark}

\subsection{Link with holes of a maximal partial spread}\label{sec:LinkHolesMaxPartSpread}
As mentioned in Section \ref{sec:constructionMaxPartSpread}, sets of points and planes fulfilling Property \ref{prop}, also appear when studying the set of holes of a maximal partial spread in $\PG(3,q)$ of size $q^2+1-x$.


 Via the study of these corresponding sets of holes, lower bounds on the size of positive values of $x$, for which a maximal partial spread of $\PG(3,q)$ of size $q^2+1-x$ exists, were found, together with many characterisation results on these sets of holes of a maximal partial spread of $\PG(3,q)$, of size $q^2+1-x$, for small $x$. 
 These results, which we will mention later on, rely on the theory of blocking sets in the plane $\PG(2,q)$.  Furthermore, the proofs for these results only relied on Property \ref{prop}. This implies that these characterisation results can be interpreted in the context of Cameron-Liebler sets of generators on the Klein quadric.

We first state them as a result on the set of holes of a maximal partial spread of size $q^2+1-x$, and then afterwards, we adapt the result to a characterisation result on Cameron-Liebler sets of generators on the Klein quadric.

We start with some required results on blocking sets in $\PG(2,q)$.

\begin{definition}
A {\em blocking set} $B$ in $\PG(2,q)$ is a set of points, intersecting every line in at least one point.
\end{definition}

Every set of points $B$ in $\PG(2,q)$, containing a line, is a blocking set of $\PG(2,q)$. This leads to the concept of non-trivial blocking sets and minimal blocking sets.

\begin{definition} A {\em non-trivial} blocking set $B$ in $\PG(2,q)$ is a blocking set, not containing a line of $\PG(2,q)$.

A {\em minimal} blocking set $B$ in $\PG(2,q)$ is a blocking set, such that for every point $P$ of $B$, the set $B\setminus \{P\}$ is no longer a blocking set.

A {\em small} blocking set $B$ in $\PG(2,q)$ is a blocking set containing less than $3(q+1)/2$ points.
\end{definition}

Classical examples of non-trivial minimal blocking sets in $\PG(2,q)$ include:

\begin{enumerate}
\item In $\PG(2,q^2)$, a Baer subplane $\PG(2,q)$ is a minimal blocking set of size $q^2+q+1$, intersecting every line in 1 or $q+1$ points. 
\item Consider a subgeometry $\PG(3,q)$ in $\PG(3,q^3)$. Project this subgeometry from a point $P$, not belonging to the subgeometry, onto a plane skew to $P$. Then there are two possibilities for this projected subgeometry.

Either the point $P$ belongs to a line of the subgeometry $\PG(3,q)$, extended to the field of order $q^3$, and then the projected subgeometry is a blocking set of size $q^3+q^2+1$.
Or the the point $P$ does not belong to a line of the subgeometry $\PG(3,q)$, extended to the field of order $q^3$, and then the projected subgeometry is a blocking set of size $q^3+q^2+q+1$.
\item These preceding examples are specific examples of linear blocking sets. The linearity conjecture states that every minimal blocking set $B$ in $\PG(2,q)$, of size $|B|< 3(q+1)/2$, is a linear blocking set, see \cite{Sziklaiblock}.

\end{enumerate}

The following results are known on this linearity conjecture and on small minimal blocking sets.

\begin{theorem}[{Sz\H onyi \cite{MR1459823}, Sziklai \cite{Sziklaiblock}}]
  Let $B$ be a minimal blocking set in $\PG(2,q), q=p^h, p$ prime, $h\geq 1$, of size $|B|< \frac{3(q+1)}{2}$. Then $B$ intersects every line in $1 \pmod{p^e}$ points. Let $e$ be the maximal integer for which this is valid. Then  $e$ divides $h$ \cite{Sziklaiblock} and  $|B|\geq q+1+\frac{q}{p^e+2}$ \cite{MR1459823}. 
\end{theorem}

\begin{theorem} 
\begin{enumerate}[label=\alph*)]
    \item {\em (Blokhuis \cite{MR1273203})}  The only small minimal blocking sets in $\PG(2,q)$, $q>2$ prime, are the lines.
    \item {\em (Sz\H onyi \cite{MR1459823})} Every small minimal non-trivial blocking set in $\PG(2,q^2)$, $q$ prime, is equal to a Baer subplane.
    \item {\em (Polverino and Storme \cite{MR1725553,MR1779313,MR1878779})} Every small minimal non-trivial blocking set in $\PG(2,q^3)$, $q\geq 7$ prime, is equal to a projected subgeometry $\PG(3,q)$ of size $q^3+q^2+1$ or of size $q^3+q^2+q+1$, or to a Baer subplane if $q$ is square.
    \item {\em (Sz\H onyi \cite{MR1459823}, Sziklai \cite{Sziklaiblock})}  Let $B$ be a minimal blocking set in $\PG(2,q), q=p^h, p$ prime, $h\geq 1$, of size $|B|< \frac{3(q+1)}{2}$, not of type $a), b)$ or $c)$. Then $|B|\geq q+1+\frac{q}{p^e+2}$, where $e$ is the largest divisor of $h$ smaller than $h/3$.
\end{enumerate}
\end{theorem}

These results on blocking sets were used to characterise the set of holes of maximal partial spreads.
We first present results in $\PG(3,q)$, $q$ square.

\begin{theorem} {\em (Blokhuis and Metsch \cite{BlokhuisMetsch}, Metsch and Storme \cite{MetschStorme})}\label{holes1}
Let $q+\delta$ be the size of the smallest non-trivial blocking set in $\PG(2,q)$, $q$ square, not containing a Baer subplane, then $\delta\leq \frac{q}{2}+1$ if $q$ is even and $\delta\leq \frac{q+1}{2}+1$ if $q$ is odd.

If $S$ is a maximal partial spread of size $q^2+1-x$, with $x<\delta$, then the set of holes of $S$ is a set of points equal to the union of $x$ Baer subgeometries $\PG(3, \sqrt{q})$, which are pairwise disjoint.

If $q>4$, then $x\geq 2(\sqrt{q}+1)$, and if $q=4$, then $x=\sqrt{q}+1=3$ is possible.
\end{theorem}

Let $q+\delta$ be the size of the smallest non-trivial blocking set in $\PG(2,q)$, $q$ square, not containing a Baer subplane, then $\delta\leq \frac{q}{2}+1$ if $q$ is even and $\delta\leq \frac{q+1}{2}+1$ if $q$ is odd, see \cite[Lemma 13.6]{Hirschfeld1}.

The 
corresponding characterisation result, {derived from Theorem \ref{holes1},} on Cameron-Liebler sets of generators in the Klein quadric is  as follows.

\begin{theorem}\label{CL1}
Let $q+\delta$ be the size of the smallest non-trivial blocking set in $\PG(2,q)$, $q$ square, not containing a Baer subplane.
Let $\L$ be a Cameron-Liebler set of generators on the Klein quadric, with $x<\delta$, then $\L$ is the union of point-pencils and of Cameron-Liebler sets of Baer subgeometry type, which are pairwise disjoint.

\end{theorem}

\begin{proof}
We want to apply {the results of} Theorem \ref{holes1} in the context of Cameron-Liebler sets. A set of holes of a maximal partial spread cannot contain a line, due to the maximality. {First consider the union of all point-pencils derived from the holes, next}  exclude a line from the set of holes of a partial spread{. This} corresponds to excluding a point-pencil as part of a Cameron-Liebler set.

To use the same proofs as for Theorem \ref{holes1}, we first assume that there is a line $\ell$ contained in the set $P_0$. Since $x<q$, this implies that all planes of $\PG(3,q)$ through this line $\ell$ belong to the set of planes $P_2$. Hence, this set of points on the line $\ell$ and this set of planes through the line $\ell$ define a Cameron-Liebler set of generators, with parameter 1, on the Klein quadric. This specific Cameron-Liebler set, with parameter 1, can be removed from the Cameron-Liebler set $\L$ with parameter $x$, to remain with a Cameron-Liebler set $\L'$ with parameter $x-1$.

This procedure can be repeated until the set of points $P_0$ does not contain a line of $\PG(3,q)$. Then the set of points $P_0$ and the set of planes $P_2$, satisfying Property \ref{prop}, has the same set of properties as the set of holes of a maximal partial spread in $\PG(3,q)$. The proofs of Theorem \ref{holes1} can be repeated to obtain Theorem \ref{CL1}.
\end{proof}\\

We now present results by S. Ferret and L. Storme. Here we use the following notation \cite{FerretStorme}.
\begin{itemize}
    \item For $q=p^3$, $p$ prime, $p\geq 17$, let 
    $\delta_0 = \lfloor (3p^3+27p^2-5p+25)/25 \rfloor$.
    \item For $p=7,11,13$, $\delta_0=90, \delta_0=285$ and $\delta_0=441$ respectively.
    \item 
For $q=p^3$, $p=p_0^h$, $p_0$ prime, $p_0 \geq 7$, $h>1$, 
let $\delta_0 = \min\{\lfloor (3p^3+27p^2-5p+25)/25 \rfloor, \delta'\}$  for which $p^3+\delta'$ is the cardinality  of the smallest non-trivial minimal blocking set in $\PG(2,p^3)$ of cardinality larger than $p^3+p^2+p+1$. Presently, this value is still unknown, but we know $\delta' \leq p^3/p_0+1$.
\end{itemize}

\begin{theorem}\label{final1} {\em (Ferret and Storme \cite{FerretStorme})}
Let $p=p_0^h$, $p_0 \geq 7$ a prime, $h\geq1$ odd.
Then the set of holes corresponding to a maximal partial spread
in $\PG(3,p^3)$ of deficiency $\delta \leq \delta_0$,
is the  union of pairwise disjoint projected subgeometries $\PG(5,p)$ of cardinality $p^5+p^4+p^3+p^2+p+1$.
\end{theorem}

\begin{theorem}\label{final2} {\em (Ferret and Storme \cite{FerretStorme})}
Let $p=p_0^h$, $p_0 \geq 7$ a prime, $h>1$ even.
Then the set of holes corresponding to a maximal partial spread
in $\PG(3,p^3)$ of deficiency $\delta \leq \delta_0$,
is the  union of pairwise disjoint Baer subgeometries $\PG(3,p^{3/2})$ and of projected subgeometries $PG(5,p)$ of cardinality $p^5+p^4+p^3+p^2+p+1$.
\end{theorem}

{The result on Cameron-Liebler sets derived from Theorem \ref{final1} and Theorem \ref{final2} is as follows.}

\begin{theorem}\label{CL2}
\begin{enumerate}
    \item Let $p=p_0^h$, $p_0 \geq 7$ a prime, $h\geq1$ odd.

Let $\L$ be a Cameron-Liebler set of generators on the Klein quadric $Q^+(5,p ^3)$, with $x\leq\delta_0$, then $\L$ is the union of point-pencils and of Cameron-Liebler sets of projected $\PG(5,p)$  type, which are pairwise disjoint.
    \item Let $p=p_0^h$, $p_0 \geq 7$ a prime, $h>1$ even.

Let $\L$ be a Cameron-Liebler set of generators on the Klein quadric $Q^+(5,p ^3)$, with $x\leq\delta_0$, then $\L$ is the union of point-pencils, Cameron-Liebler sets of Baer subgeometry type,  and of Cameron-Liebler sets of projected $\PG(5,q)$  type, which are pairwise disjoint.
\end{enumerate}

\end{theorem}

\bibliographystyle{abbrv}
\bibliography{sample.bib}

\begin{thebibliography}{10}

\bibitem{MR1273203}
A.~Blokhuis.
\newblock On the size of a blocking set in {${\rm PG}(2,p)$}.
\newblock {\em Combinatorica}, 14(1):111--114, 1994.

\bibitem{bdd}
A.~Blokhuis, M.~{De Boeck}, and J.~D'haeseleer.
\newblock {Cameron--Liebler} sets of $k$-spaces in {$\PG(n,q)$}.
\newblock {\em Des. Codes Cryptogr.}, 87:1839--1856, 2018.

\bibitem{bddc}
A.~Blokhuis, M.~{De Boeck}, and J.~D'haeseleer.
\newblock Correction to: {Cameron--Liebler} sets of $k$-spaces in {$\PG(n,q)$}.
\newblock {\em Des. Codes Cryptogr.}, 90:477--487, 2022.

\bibitem{MR1772206}
A.~Blokhuis and M.~Lavrauw.
\newblock Scattered spaces with respect to a spread in {${\rm PG}(n,q)$}.
\newblock {\em Geom. Dedicata}, 81(1-3):231--243, 2000.

\bibitem{BlokhuisMetsch}
A.~Blokhuis and K.~Metsch.
\newblock On the size of a maximal partial spread.
\newblock {\em Des. Codes Cryptogr.}, 3(3):187--191, 1993.

\bibitem{BCN}
A.~E. Brouwer, A.~M. Cohen, and A.~Neumaier.
\newblock {\em Distance-regular graphs}, volume~18 of {\em Ergebnisse der Mathematik und ihrer Grenzgebiete (3) [Results in Mathematics and Related Areas (3)]}.
\newblock Springer-Verlag, Berlin, 1989.

\bibitem{bd}
A.~Bruen and K.~Drudge.
\newblock The construction of {C}ameron-{L}iebler line classes in {$\PG(3,q)$}.
\newblock {\em Finite Fields Appl.}, 5:35--45, 1999.

\bibitem{cameronliebler}
P.~Cameron and R.~Liebler.
\newblock Tactical decompositions and orbits of projective groups.
\newblock {\em Linear Algebra Appl.}, 46:91--102, 1982.

\bibitem{cossidentepavese}
A.~Cossidente and F.~Pavese.
\newblock {Cameron--Liebler} line classes of {$\PG(3,q)$} admitting {$\PGL(2,q)$}.
\newblock {\em J. Combin. Theory Ser. A}, 167:104--120, 2019.

\bibitem{DBDIM23}
J.~De~Beule, J.~D'haeseleer, F.~Ihringer, and J.~Mannaert.
\newblock Degree 2 {B}oolean functions on {G}rassmann graphs.
\newblock {\em Electron. J. Combin.}, 30(1):Paper No. 1.31, 23, 2023.

\bibitem{debeulemannaert}
J.~{De Beule} and J.~Mannaert.
\newblock A modular equality for {Cameron-Liebler} line classes in projective and affine spaces of odd dimension.
\newblock {\em Finite Fields Appl.}, 82:102047, 2022.

\bibitem{debeulemannaertstorme}
J.~De~Beule, J.~Mannaert, and L.~Storme.
\newblock {Cameron--Liebler} $k$-sets in subspaces and non-existence conditions.
\newblock {\em Des. Codes Cryptogr.}, 90:633--651, 2022.

\bibitem{CLpolequiv}
M.~De~Boeck and J.~D'haeseleer.
\newblock Equivalent definitions for (degree one) {C}ameron-{L}iebler classes of generators in finite classical polar spaces.
\newblock {\em Discrete Math.}, 343(1):111642, 13, 2020.

\bibitem{CLpolMMJ}
M.~De~Boeck, J.~D'haeseleer, and M.~Rodgers.
\newblock Regular ovoids and {C}ameron-{L}iebler sets of generators in polar spaces.
\newblock {\em J. Combin. Theory Ser. A}, 213:Paper No. 106029, 2025.

\bibitem{CLpol}
M.~De~Boeck, M.~Rodgers, L.~Storme, and A.~\v{S}vob.
\newblock Cameron-{L}iebler sets of generators in finite classical polar spaces.
\newblock {\em J. Combin. Theory Ser. A}, 167:340--388, 2019.

\bibitem{debruynsuzuki}
B.~{De Bruyn} and H.~Suzuki.
\newblock Intriguing sets of vertices of regular graphs.
\newblock {\em Graphs Combin.}, 26:629--646, 2010.

\bibitem{D99}
K.~Drudge.
\newblock On a {C}onjecture of {C}ameron and {L}iebler.
\newblock {\em European J. Combin}, 20:263--269, May 1999.

\bibitem{ThesisDrudge}
K.~W. Drudge.
\newblock {\em Extremal sets in projective and polar spaces}.
\newblock ProQuest LLC, Ann Arbor, MI, 1998.
\newblock Thesis (Ph.D.)--The University of Western Ontario (Canada).

\bibitem{TF:14}
T.~Feng, K.~Momihara, and Q.~Xiang.
\newblock {Cameron-Liebler} line classes with parameter $x=\frac{q^2-1}{2}$.
\newblock {\em J. Combin. Theory Ser. A}, 133:307--338, 2015.

\bibitem{FerretStorme}
S.~Ferret and L.~Storme.
\newblock Results on maximal partial spreads in {${\rm PG}(3,p^3)$} and on related minihypers.
\newblock {\em Des. Codes Cryptogr.}, 29(1-3):105--122, 2003.

\bibitem{Ferdinand.}
Y.~Filmus and F.~Ihringer.
\newblock Boolean degree 1 functions on some classical association schemes.
\newblock {\em J. Combin. Theory Ser. A}, 162:241--270, 2019.

\bibitem{gmp}
A.~Gavrilyuk, I.~Matkin, and T.~Penttila.
\newblock Derivation of {Cameron--Liebler} line classes.
\newblock {\em Des. Codes Cryptogr.}, 86:231--236, 2018.

\bibitem{gmet}
A.~Gavrilyuk and K.~Metsch.
\newblock A modular equality for {Cameron-Liebler} line classes.
\newblock {\em J. Combin. Theory Ser. A}, 127:224--242, 2014.

\bibitem{gmol}
A.~Gavrilyuk and I.~Mogilnykh.
\newblock {Cameron-Liebler} line classes in {$\PG(n,4)$}.
\newblock {\em Des. Codes Cryptogr.}, 73:969--982, 2014.

\bibitem{CLaffineclassical}
J.~Guo and L.~Wan.
\newblock Cameron-{L}iebler sets for maximal totally isotropic flats in classical affine spaces.
\newblock {\em J. Combin. Des.}, 31(11):547--574, 2023.

\bibitem{MR1235896}
O.~Heden.
\newblock Maximal partial spreads and the modular {$n$}-queen problem.
\newblock {\em Discrete Math.}, 120(1-3):75--91, 1993.

\bibitem{MR1341439}
O.~Heden.
\newblock Maximal partial spreads and the modular {$n$}-queen problem. {II}.
\newblock {\em Discrete Math.}, 142(1-3):97--106, 1995.

\bibitem{MR1874731}
O.~Heden.
\newblock Maximal partial spreads and the modular {$n$}-queen problem. {III}.
\newblock {\em Discrete Math.}, 243(1-3):135--150, 2002.

\bibitem{HedenStorme}
O.~Heden and L.~Storme.
\newblock Maximal partial spreads and the modular {$n$}-queen problem for $q$ even (unpublished manuscript).

\bibitem{Hirschfeld1}
J.~W.~P. Hirschfeld.
\newblock {\em Projective geometries over finite fields}.
\newblock Oxford Mathematical Monographs. The Clarendon Press, Oxford University Press, New York, second edition, 1998.

\bibitem{PhDLavrauw}
M.~Lavrauw.
\newblock Scattered spaces with respect to spreads and eggs in finite projective spaces, phd thesis, {E}indhoven {U}niversity of {T}echnology, 2001.

\bibitem{met}
K.~Metsch.
\newblock The non-existence of {Cameron-Liebler} line classes with parameter $2<x\leq q$.
\newblock {\em Bull. Lond. Math. Soc.}, 42:991--996, 2010.

\bibitem{met2}
K.~Metsch.
\newblock An improved bound on the existence of {Cameron-Liebler} line classes.
\newblock {\em J. Combin. Theory Ser. A}, 121:89--93, 2014.

\bibitem{met3}
K.~Metsch.
\newblock A gap result for {Cameron-Liebler} $k$-classes.
\newblock {\em Discrete Math.}, 340:1311--1318, 2017.

\bibitem{MetschStorme}
K.~Metsch and L.~Storme.
\newblock Partial {$t$}-spreads in {${\rm PG}(2t+1,q)$}.
\newblock volume~18, pages 199--216. 1999.
\newblock Designs and codes---a memorial tribute to Ed Assmus.

\bibitem{meysets}
A.~Meyerowitz.
\newblock Cycle-balanced partitions in distance-regular graphs.
\newblock {\em J. Comb. Inf. Syst. Sci.}, 17:39--42, 1992.

\bibitem{MR1725553}
O.~Polverino.
\newblock Small minimal blocking sets and complete {$k$}-arcs in {${\rm PG}(2,p^3)$}.
\newblock {\em Discrete Math.}, 208/209:469--476, 1999.

\bibitem{MR1779313}
O.~Polverino.
\newblock Small blocking sets in {${\rm PG}(2,p^3)$}.
\newblock {\em Des. Codes Cryptogr.}, 20(3):319--324, 2000.

\bibitem{MR1878779}
O.~Polverino and L.~Storme.
\newblock Small minimal blocking sets in {${\rm PG}(2,q^3)$}.
\newblock {\em European J. Combin.}, 23(1):83--92, 2002.

\bibitem{rod}
M.~Rodgers.
\newblock {Cameron-Liebler} line classes.
\newblock {\em Des. Codes Cryptogr.}, 68:33--37, 2013.

\bibitem{rodstovan}
M.~Rodgers, L.~Storme, and A.~Vansweevelt.
\newblock {Cameron-Liebler} $k$-classes in {$\PG(2k+1,q)$}.
\newblock {\em Combinatorica}, 38:739--757, 2018.

\bibitem{MR1459823}
T.~Sz\H{o}nyi.
\newblock Blocking sets in {D}esarguesian affine and projective planes.
\newblock {\em Finite Fields Appl.}, 3(3):187--202, 1997.

\bibitem{Sziklaiblock}
P.~Sziklai.
\newblock On small blocking sets and their linearity.
\newblock {\em J. Combin. Theory Ser. A}, 115(7):1167--1182, 2008.

\bibitem{tits}
J.~Tits.
\newblock {\em Buildings of spherical type and finite {BN}-pairs}, volume 386 of {\em Lecture Notes in Mathematics}.
\newblock Springer-Verlag, Berlin-New York, 1974.

\bibitem{Veldkamp}
F.~D. Veldkamp.
\newblock Polar geometry, {I--V}.
\newblock {\em Proc. Kon. Ned. Akad. Wet.}, A62:512--551, 1959.

\end{thebibliography}
\newpage

\appendix

\end{document}